\documentstyle[twoside]{article}
\include{graphicx}
\title{\bf Asymptotics of a Gauss hypergeometric function with large parameters, III: Application to the Legendre functions of large imaginary order and real degree}
\author{\sc R. B. Paris\footnote{E-mail address:\ \ {\tt r.paris@abertay.ac.uk}}\\
\\
{\em Division of Computing and Mathematics,}\\
{\em University of Abertay Dundee, Dundee DD1 1HG, UK}\\
}

\begin{document}
\newcommand{\bee}{\begin{equation}}
\newcommand{\ee}{\end{equation}}
\def\f#1#2{\mbox{${\textstyle \frac{#1}{#2}}$}}
\def\dfrac#1#2{\displaystyle{\frac{#1}{#2}}}
\newcommand{\fr}{\frac{1}{2}}
\newcommand{\fs}{\f{1}{2}}
\newcommand{\g}{\Gamma}
\newcommand{\br}{\biggr}
\newcommand{\bl}{\biggl}
\newcommand{\ra}{\rightarrow}
\renewcommand{\topfraction}{0.9}
\renewcommand{\bottomfraction}{0.9}
\renewcommand{\textfraction}{0.05}
\newcommand{\mcol}{\multicolumn}
\date{}
\maketitle
\pagestyle{myheadings}
\markboth{\hfill {\it R.B. Paris} \hfill}
{\hfill {\it Hypergeometric function with large parameters} \hfill}
\begin{abstract} 
We obtain the asymptotic expansion for the Gauss hypergeometric function
\[F(a-\lambda,b+\lambda;c+i\alpha\lambda;z)\]
for $\lambda\ra+\infty$ with $a$, $b$ and $c$ finite parameters by application of the method of steepest descents. The quantity $\alpha$ is real, so that the denominatorial parameter is complex and $z$ is a finite complex variable restricted to lie in the sector $|\arg (1-z)|<\pi$. We concentrate on the particular case $a=0$, $b=c=1$, which is associated with the Legendre functions of real degree and imaginary order. The resulting expansions are of Poincar\'e type and hold in restricted domains of the $z$-plane. An expansion is given at the coalescence of two saddle points. Numerical results illustrating the accuracy of the different expansions are given. 
\vspace{0.4cm}

\noindent {\bf MSC:} 33C05, 34E05, 41A60
\vspace{0.3cm}

\noindent {\bf Keywords:} Hypergeometric function, asymptotic expansion, large parameters, steepest descents, Legendre functions\\
\end{abstract}

\vspace{0.2cm}

\noindent $\,$\hrulefill $\,$

\vspace{0.2cm}

\begin{center}
{\bf 1. \  Introduction}
\end{center}
\setcounter{section}{1}
\setcounter{equation}{0}
\renewcommand{\theequation}{\arabic{section}.\arabic{equation}}
The Gauss hypergeometric function is defined by
\[F\left(\!\!\begin{array}{c}a, b\\c\end{array}\!;z\!\right)=\sum_{n=0}^\infty\frac{(a)_n (b)_n}{(c)_n n!}\,z^n \qquad(|z|<1)\]
and elsewhere by analytic continuation, where $(a)_n=\g(a+n)/\g(a)=a(a+1)\ldots (a+n-1)$ is the Pochhammer symbol or rising factorial. The asymptotic expansion of 
\[F\left(\begin{array}{c}a+\epsilon_1\lambda, b+\epsilon_2\lambda\\c+\epsilon_3\lambda\end{array}\!;z\right)\]
for large values of $\lambda$ and fixed complex $z$ when the parameters $\epsilon_j$ are finite was first considered by Watson \cite{W} in 1918  and recently by the author in \cite{P1, P2}; see also \cite{Cvit} for the case of two large parameters.  This function may be characterised by the set $\{\epsilon_1, \epsilon_2, \epsilon_3\}$, where, by a rescaling of $\lambda$ one of the $\epsilon_j$ can (if so desired) always be replaced by unity. Watson considered the situation when $\epsilon_j=0, \pm 1$ and examined the cases $(0,0,1)$, $(1,-1,0)$ and $(0,-1,1)$, together with the additional case $(1,2,0)$.
In \cite{P1, P2}, the $\epsilon_j$ were taken to be positive constants. It was shown that it is sufficient to consider just three basic types of hypergeometric function corresponding to the parameter sets $(\epsilon,0,1)$ (Type A), $(\epsilon,-1,0)$ (Type B) and $(\epsilon_1, \epsilon_2,1)$ (Type C), where $\epsilon, \epsilon_1, \epsilon_2>0$. An application of the expansion for the case
$(\epsilon, \epsilon,1)$ with $\epsilon>0$ has arisen in aerodynamics \cite{Cherry, Lig}. Expansions of a uniform character when two parameters are large have been given for the case $(1, -1,0)$ in \cite{DSJ} and $(0, -1, 1)$, $(1, 2, 0)$ in \cite{OD}. 

In this paper, we consider the expansion of the function
\bee\label{e11}
F\left(\begin{array}{c}a-\lambda, b+\lambda\\c+i\alpha\lambda\end{array}\!;z\right)
\ee
for $\lambda\ra+\infty$ with finite $a$, $b$ and $c$. 
Here the parameter $\alpha\in {\bf R}$ is finite and $z$ is a (finite) complex variable restricted to lie in $|\arg (1-z)|<\pi$.
The parameter set in this case is consequently $(-1,1,i\alpha)$, which is different from the cases previously considered in that the denominatorial parameter $\epsilon_3=i\alpha$ is {\it purely imaginary}. 
As in \cite{P1,P2}, we employ the method of steepest descents applied to a contour integral representation of the function in (\ref{e11}) to obtain Poincar\'e-type expansions. The expansion in the case of the coalescence of two saddle points is also considered.

We concentrate on the situation when the parameters in (\ref{e11})
have the values $a=0$, $b=c=1$ and consider the function
\bee\label{e20}
F_\alpha(\lambda;z):=F\left(\begin{array}{c}-\lambda, 1+\lambda\\1+i\alpha\lambda\end{array}\!;z\right)
\ee
with $\alpha>0$.
This particular case has arisen in the study of travelling waves in a Toda lattice \cite{Wr,VSWP}. 
This case is also of interest as it is associated with the Legendre functions through the relations \cite[p.~353]{DLMF}
\bee\label{e12}
P_\nu^{-\mu}(x)=\frac{1}{\g(1+\mu)} \left(\frac{x-1}{x+1}\right)^{\mu/2} F(-\nu,1+\nu;1+\mu;\fs-\fs x)
\ee
and 
\[e^{\pi i\mu}Q_\nu^{-\mu}(x)=\frac{\g(1+\nu-\mu)\g(\mu)}{2\g(1+\nu+\mu)}\left(\frac{x-1}{x+1}\right)^{-\mu/2} F(-\nu,1+\nu;1-\mu;\fs-\fs x)\]
\bee\label{e13}
\hspace{4cm}+\frac{\g(-\mu)}{2} \left(\frac{x-1}{x+1}\right)^{\mu/2} F(-\nu,1+\nu;1+\mu;\fs-\fs x),
\ee
which define the functions in the complex $x$-plane cut along $(-\infty,1]$.
Thus, as a by-product of our investigation of (\ref{e20}) we will obtain the expansions of the Legendre functions $P_\lambda^{\pm i\alpha\lambda}(x)$ and $Q_\lambda^{\pm i\alpha\lambda}(x)$, for large imaginary order and real degree.
The expansions for these functions when $x=\sqrt{1+\alpha^2}$, which corresponds to the above-mentioned coalescence of two saddle points, is also given. We remark that the expansion of $P_\nu^{-i\mu}(x)$ and $Q_\nu^{-i\mu}(x)$ for $\nu\ra+\infty$, $\Re (x)\geq 0$ uniformly valid for $0\leq \mu/(\nu+\fs)\leq A$, where $A$ is a constant, has been considered by Dunster \cite{TMD} who employed a differential-equation approach.

We remark that when $z=\fs$ it is possible to give an exact evaluation for $F_\alpha(\lambda;z)$ in the form \cite[(15.4.30)]{DLMF}
\bee\label{e14}
F_\alpha(\lambda;\fs)=\frac{2^{-i\alpha\lambda}\sqrt{\pi}\,\g(1+i\alpha\lambda)}{\g(\fs+\fs\lambda(i\alpha-1))\g(1+\fs\lambda(i\alpha+1))}~.
\ee
\vspace{0.6cm}

\begin{center}
{\bf 2. \ The expansion of $F_\alpha(\lambda;z)$ for $\lambda\ra+\infty$}
\end{center}
\setcounter{section}{2}
\setcounter{equation}{0}
\renewcommand{\theequation}{\arabic{section}.\arabic{equation}}
We take the parameter values  $a=0$, $b=c=1$ in (\ref{e11})   and consider in detail the expansion of $F_\alpha(\lambda;z)$ in (\ref{e20}) for $\lambda\ra+\infty$; 
the case of general values of $a$, $b$ and $c$ is considered in Section 2.3.
From the series representation of the hypergeometric function it follows that, when $\lambda>0$ and $\alpha$ is real,
\bee\label{e21a}
F_{-\alpha}(\lambda;z)={\overline F}_\alpha(\lambda;{\overline z}),
\ee
where the bar denotes the complex conjugate. It is therefore sufficient to consider $\alpha>0$ throughout; in addition, we define
\bee\label{e200}
\theta:=\arg\,z,\qquad \phi:=\arctan\,\alpha,
\ee
where $\theta\in [-\pi,\pi]$ and $\phi\in (0,\fs\pi)$.

We employ the integral representation \cite[p.~388]{DLMF}
\bee\label{e2int}
F\left(\!\!\begin{array}{c}a, b\\c\end{array}\!;z\!\right)=
\frac{\g(c) \g(1+b-c)}{2\pi i \g(b)}\int_0^{(1+)} \frac{t^{b-1}(t-1)^{c-b-1}}{(1-zt)^a}\,dt, \qquad \Re (b)>0,
\ee
where it is supposed that $|\arg (1-z)|<\pi$ and $c-b\neq 1, 2, \ldots\ $. The integration path is a closed loop that starts from the origin, encircles the point $t=1$ in the positive sense and returns to the origin without 
enclosing the point $t=1/z$. 
The function $F_\alpha(\lambda;z)$ can then be expressed in the form
\bee\label{e21}
F_\alpha(\lambda;z)=F\left(\begin{array}{c}-\lambda, 1+\lambda\\1+i\alpha\lambda\end{array}\!;z\right)=\frac{G_\alpha(\lambda)}{2\pi i}\int_0^{(1+)} f(t) e^{\lambda\psi(t)} dt,
\ee
where the phase function $\psi(t)$ and the amplitude function $f(t)$ are
\bee\label{e22}
\psi(t)=\log\,t(1-zt)-\beta \log (t-1),\qquad f(t)=(t-1)^{-1},\qquad \beta:=1-i\alpha
\ee
and 
\bee\label{e23}
G_\alpha(\lambda):=\frac{\g(1+i\alpha\lambda) \g(1+\lambda\beta)}{\g(1+\lambda)}=i\alpha\beta\lambda\,\frac{\g(i\alpha\lambda)\g(\lambda\beta)}{\g(\lambda)}~.
\ee 
The $t$-plane is cut along $(-\infty,1]$ and also along the ray from the singularity at $t=1/z$ to infinity in a suitable direction.

The phase function has saddle points where $\psi'(t)=0$; that is at the points where
\[\frac{i\alpha t_s-1}{t_s(t_s-1)}-\frac{z}{1-zt_s}=0.\]
There are consequently two saddle points, which we label $t_{s1}$ and $t_{s2}$, given by
\bee\label{e24}
t_{s1}, t_{s2}=\frac{z+\fs i\alpha \mp i(z-z^2+\f{1}{4}\alpha^2)^{1/2}}{(1+i\alpha)z},
\ee
respectively. For sufficiently large $\lambda$, the points $t=0$ and $t=1/z$ are zeros of the integrand, so that paths of steepest descent can terminate only at these two points; paths of steepest ascent must terminate at $t=1$ and at infinity. A typical arrangement of the steepest paths through $t_{s1}$ and $t_{s2}$ is shown in Fig.~1 when $\alpha=1$
and for different values of $\theta$. 

Since $t_{s1}t_{s2}=\{(1+i\alpha)z\}^{-1}$, it follows that
\[\arg\,t_{s1}+\arg\,t_{s2}=-(\theta+\phi).\]
Consequently, when $\theta+\phi>0$ (resp. $<0$) at {\it least\/} one saddle is situated in the lower (resp. upper) half-plane; when $\theta+\phi=0$, one  saddle is situated in upper half-plane with the other in the lower half-plane. It is to be noted that the saddles coalesce to form a double saddle when $z^2-z-\f{1}{4}\alpha^2=0$; that is when 
\bee\label{ezd}
z=z_d^\pm:=\fs\pm \fs \sqrt{1+\alpha^2}. 
\ee
\begin{figure}[t]
	\begin{center}
{\tiny($a$)}\ \includegraphics[width=0.33\textwidth]{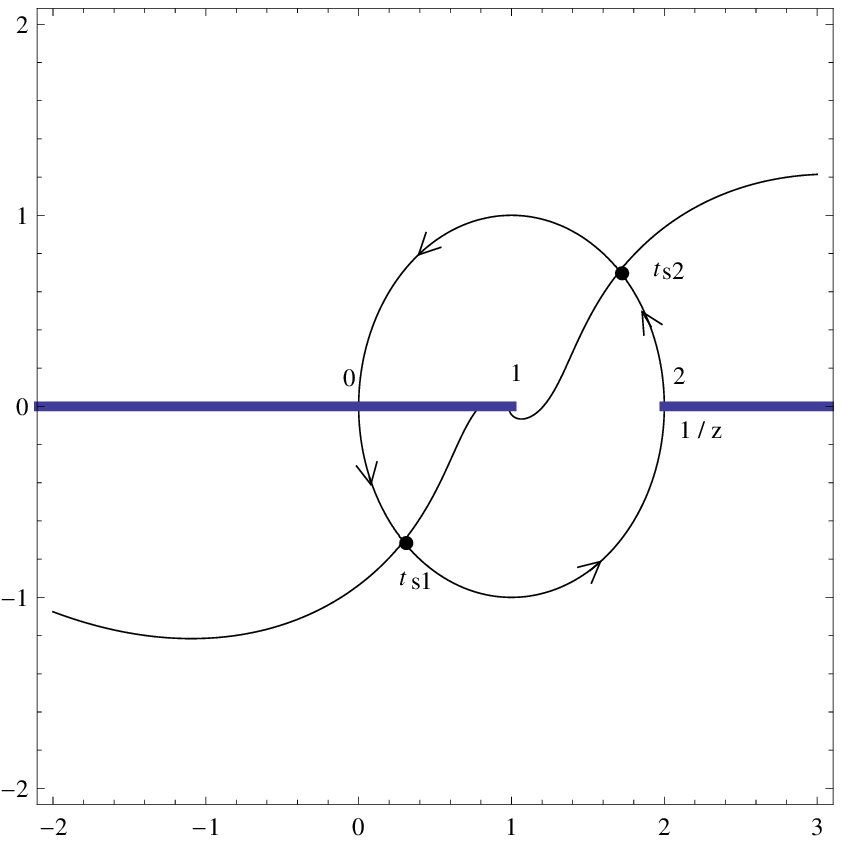}\qquad
{\tiny($b$)}\ \includegraphics[width=0.33\textwidth]{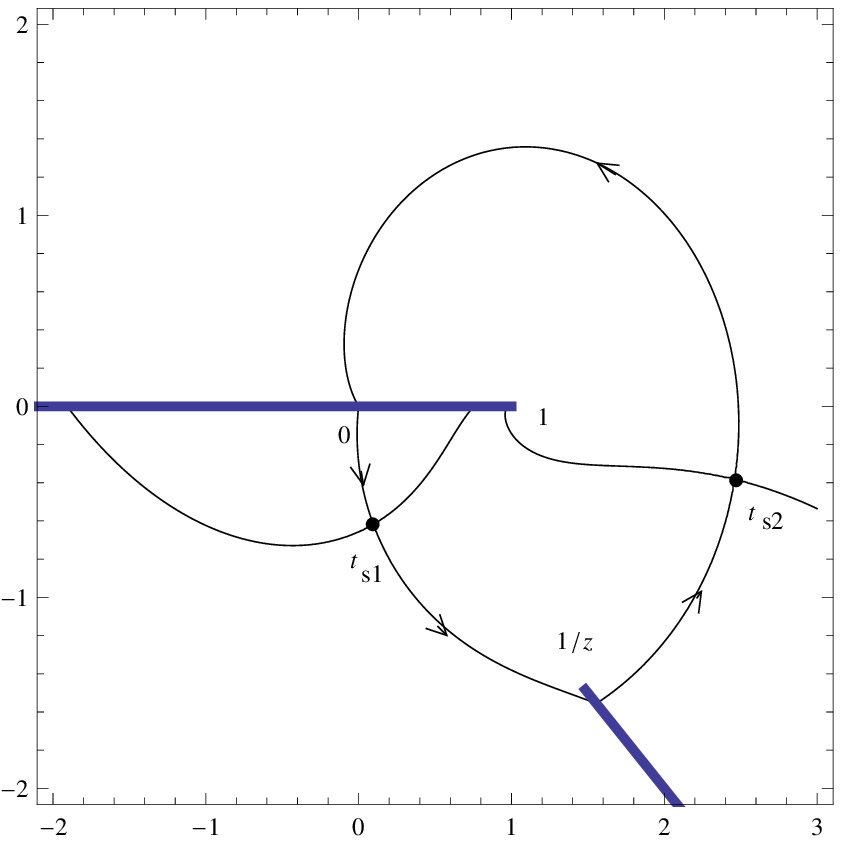}
\vspace{0.3cm}

{\tiny($c$)}\ \includegraphics[width=0.33\textwidth]{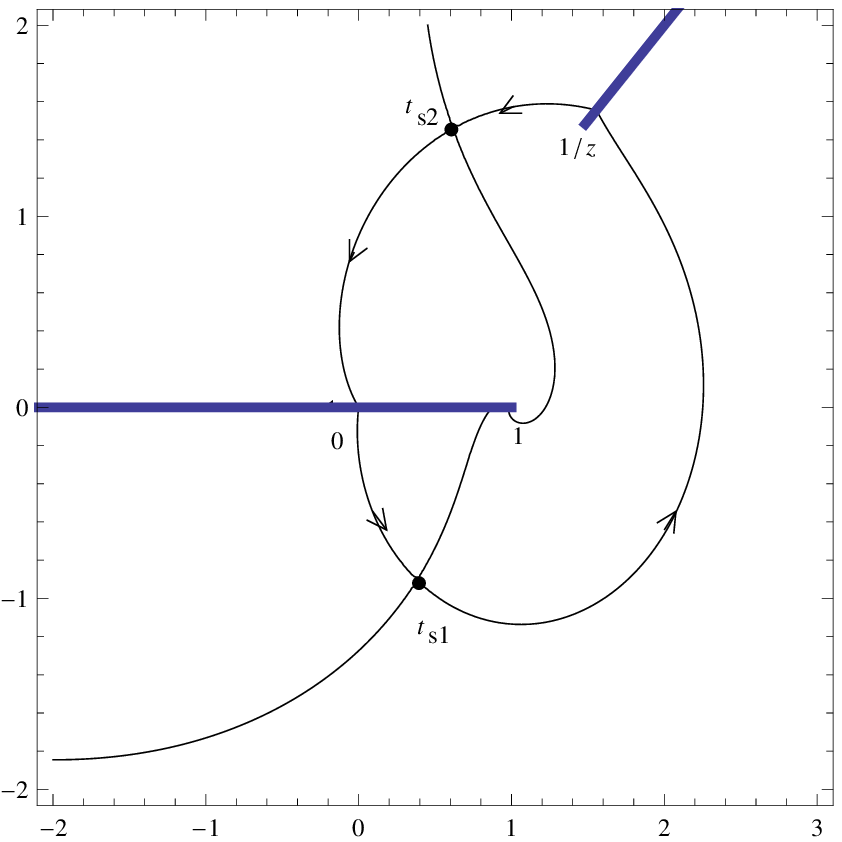}\qquad
{\tiny($d$)}\ \includegraphics[width=0.33\textwidth]{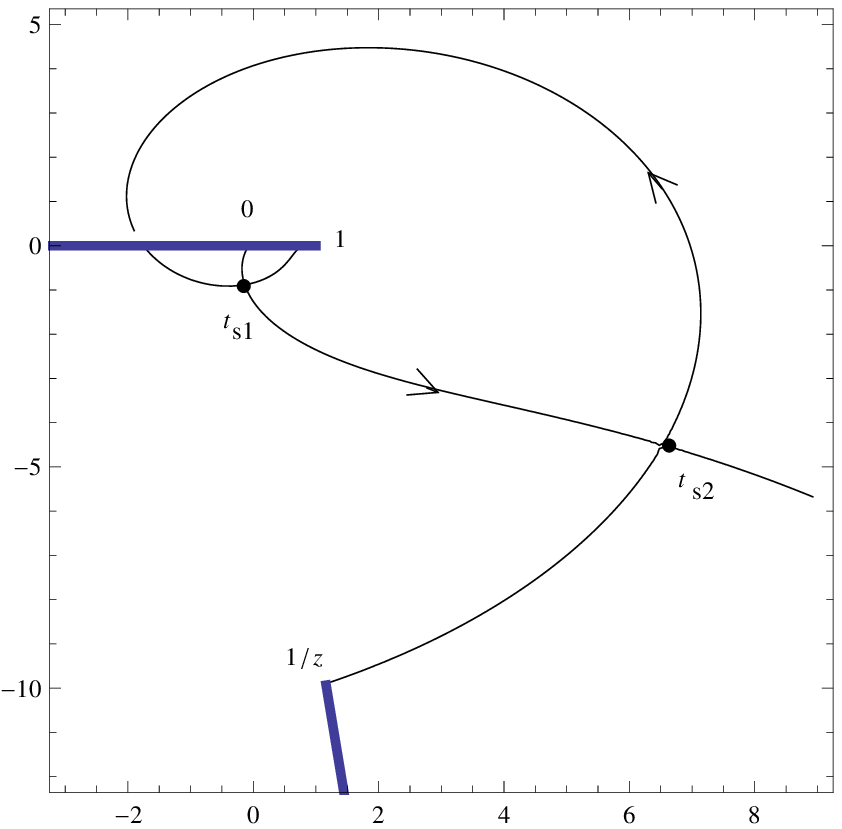}
\caption{\small{Examples of the steepest descent and ascent paths for $|z|=0.50$ and $\alpha=1$ when (a) $\theta=0$, (b) $\theta=\f{1}{4}\pi$ and (c) $\theta=-\f{1}{4}\pi$; (d) $|z|=0.10$, $\theta\doteq 0.46292\pi$ corresponding to a Stokes phenomenon. The saddles are denoted by heavy dots and the arrows indicate the integration path. There is a branch cut along $(-\infty,1]$ and from the point $1/z$ out to infinity.}}
	\end{center}
\end{figure}

The contribution to the integral (\ref{e21}) (excluding the pre-factor) from the steepest descent path through the (simple) saddle $t_{sj}$ is given by the formal asymptotic sum \cite[p.~265]{BH}
\bee\label{e25}
if(t_{sj}) e^{\lambda\psi(t_{sj})+\pi i\gamma_j}\sqrt{\frac{2\pi}{\psi''(t_{sj})}}\ \sum_{s=0}^\infty\frac{c_s^{(j)} (\fs)_s}{\lambda^{s+\fr}}\qquad (j=1, 2)
\ee
as $\lambda\ra+\infty$, where
\[\psi''(t_{sj})=-\frac{2}{(t_{sj}-1)^2}\bl\{\bl(1-\frac{1}{t_{sj}}\br)\bl(1-\frac{1}{t_{sj}}-\beta\br)-\frac{1}{2}\beta(1-\beta)\br\}.\]
The $\gamma_j$ are orientation factors that depend on the direction of integration $\arg (t-t_{sj})$ through the saddle point $t_{sj}$ and have the value either 0 or 1. 
\vspace{0.4cm}

\noindent{\bf 2.1\ \ The coefficients $c_s^{(j)}$}.\ \ \ 
The coefficients $c_s^{(j)}\equiv c_s$ (which are functions of $\alpha$ and $z$) for $s\leq 2$ are given explicitly by
\[c_0=1,\qquad c_1=\frac{1}{2\psi''}\{2F_2-2\Psi_3F_1+\f{5}{6}\Psi_3^2-\fs\Psi_4\},\]
\[c_2=\frac{1}{(2\psi'')^2}\{\f{2}{3}F_4-\f{20}{9}\Psi_3F_3+\f{5}{3}(\f{7}{3}\Psi_3^2-\Psi_4)F_2-\f{35}{9}(\Psi_3^3-\Psi_3\Psi_4+\f{6}{35}\Psi_5)F_1 \]
\bee\label{e25c}
+\f{35}{9}(\f{11}{24}\Psi_3^4-\f{3}{4}(\Psi_3^2-\f{1}{6}\Psi_4)\Psi_4+\f{1}{5}\Psi_3\Psi_5-\f{1}{35}\Psi_6)\},
\ee
where, for brevity, we have defined
\[\Psi_m:=\frac{\psi^{(m)}(t)}{\psi''(t)}\ \ \ (m\geq 3),\qquad F_m:=\frac{f^{(m)}(t)}{f(t)}\ \ \ (m\geq 1)\]
with $\psi(t)$, $f(t)$ and their derivatives being evaluated at $t=t_{s1}$ or $t=t_{s2}$; see, for example, \cite[p.~119]{D}, \cite[p.~127]{O}, \cite[p.~13]{PBook}. 

Alternatively, the $c_s$ can be obtained by an expansion process to yield Wojdylo's formula \cite{Woj} given by
\bee\label{e25w}
c_{s}=\frac{(-)^s}{{\hat\alpha}_0^{s}}\sum_{k=0}^{2s} \frac{{\hat\beta}_{2s-k}}{{\hat\beta}_0} \sum_{j=0}^k \frac{(-)^j (s+\fs)_j}{j!\ {\hat\alpha}_0^j}\,{\cal B}_{kj}\,;
\ee
see also \cite[p.~25]{T}. Here ${\cal B}_{kj}\equiv {\cal B}_{kj}({\hat\alpha}_1, {\hat\alpha}_2, \ldots , {\hat\alpha}_{k-j+1})$ are the partial ordinary Bell polynomials generated by the recursion\footnote{For example, this generates the values ${\cal B}_{41}={\hat\alpha}_4$, ${\cal B}_{42}={\hat\alpha}_2^2+2{\hat\alpha}_1{\hat\alpha}_3$, ${\cal B}_{43}=3{\hat\alpha}_1^2{\hat\alpha}_2$ and ${\cal B}_{44}={\hat\alpha}_1^4$.}
\[{\cal B}_{kj}=\sum_{r=1}^{k-j+1} {\hat\alpha}_r {\cal B}_{k-r,j-1} ,\qquad {\cal B}_{k0}=\delta_{k0},\]
where $\delta_{mn}$ is the Kronecker symbol, and the coefficients ${\hat\alpha}_r$ and ${\hat\beta}_r$ appear in the expansions
\[
\psi(t)-\psi(t_s)=\sum_{r=0}^\infty {\hat\alpha}_r (t-t_s)^{r+2},\qquad f(t)=\sum_{r=0}^\infty {\hat\beta}_r(t-t_s)^r
\]
valid in a neighbourhood of the saddle $t=t_s$.
\vspace{0.4cm}

\noindent{\bf 2.2\ \ The expansion of $F_\alpha(\lambda;z)$.}\ \ \ 
It can be seen from Fig.~1 that the {\it steepest ascent\/} path through $t_{s1}$ crosses the branch cut to pass onto an adjacent Riemann sheet in the $t$-plane. When $\theta$ increases, it is found that in some cases the saddle $t_{s2}$ (and consequently a portion of the associated steepest descent path) can also pass onto an adjacent sheet. To avoid this difficulty, we make the substitution $t=e^w$ to find the phase function in (\ref{e22}) given by
\[\psi(w)=(1-\fs\beta)w+\log (1-ze^w)-\beta \log\,2\sinh \fs w,\]
with the image of the saddles given by $w_j=\log\,t_{sj}$ ($j=1, 2$).
The closed circuit surrounding the point $t=1$ in the $t$-plane becomes the loop that commences at $-\infty$, encircles the point $w=0$ in the positive sense and returns to $-\infty$. 
Branch cuts are introduced along $(-\infty,1]$ and from the point $\log\,1/z$ out to $\infty$ parallel to the real $w$-axis. 
This transformation causes all the Riemann sheets in the $t$-plane to appear as horizontal strips of width $2\pi$ in the $w$-plane; the principal sheet corresponds to $-\pi<\Im (w)\leq\pi$. 

Examples of the steepest paths in the $w$-plane are illustrated in Fig.~2. In Fig.~2(a) both the saddles $w_1$ and $w_2$ and their associated steepest descent paths are situated on the principal sheet; in Fig.~2(b) the saddles are again on the principal sheet, but the steepest descent path from $w_2$ crosses the line $\Im (w)=\pi$, which corresponds to passing on to the adjacent sheet in the $t$-plane. Fig.~2(c) shows the same situation as Fig.~1(d), namely $|z|=0.10$, $\theta\doteq 0.46292\pi$; for this value of $\theta$ the saddles $w_1$ and $w_2$ are connected, with the steepest descent path from $w_2$ passing into the strip $\pi<\Im (w)\leq 3\pi$. The saddle in this strip corresponds to the image of the saddle $w_1$ in the principal sheet; the contribution from this image saddle is exponentially smaller (by the factor $e^{-2\pi\alpha\lambda}$) than that from $w_1$ and so is neglected. Fig.~2(d) shows $|z|=0.10$, $\theta=0.60\pi$ where the steepest descent path through $w_1$ has disconnected from the saddle $w_2$ (a Stokes phenomenon) and passes over into the adjacent strip.
\begin{figure}[t]
	\begin{center}
{\tiny($a$)}\ \includegraphics[width=0.33\textwidth]{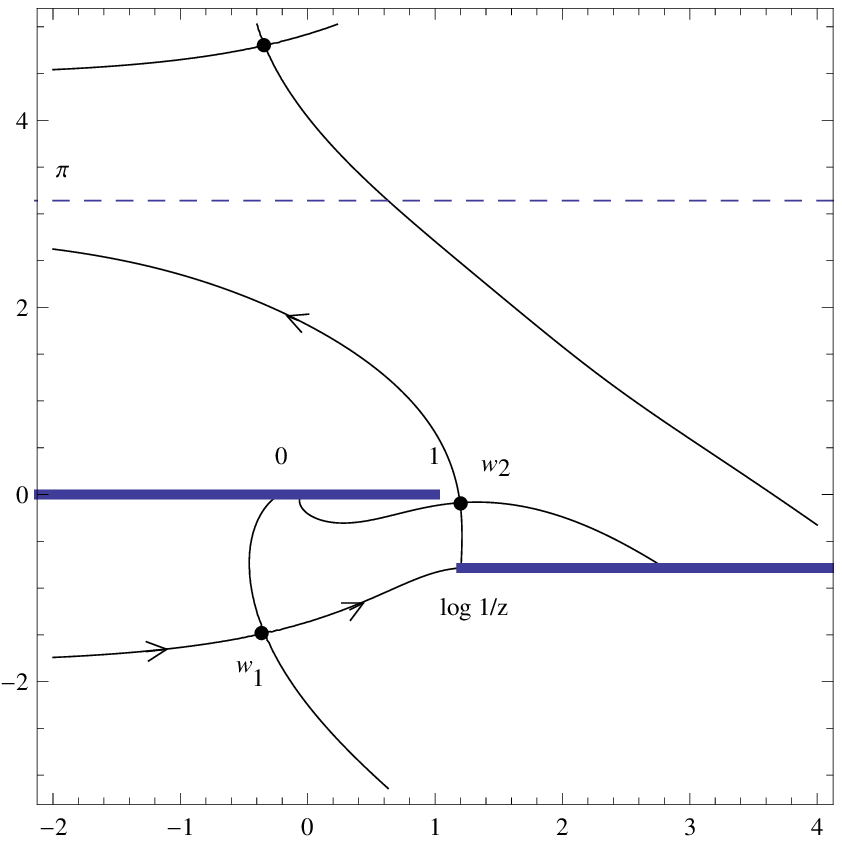}\qquad
{\tiny($b$)}\ \includegraphics[width=0.33\textwidth]{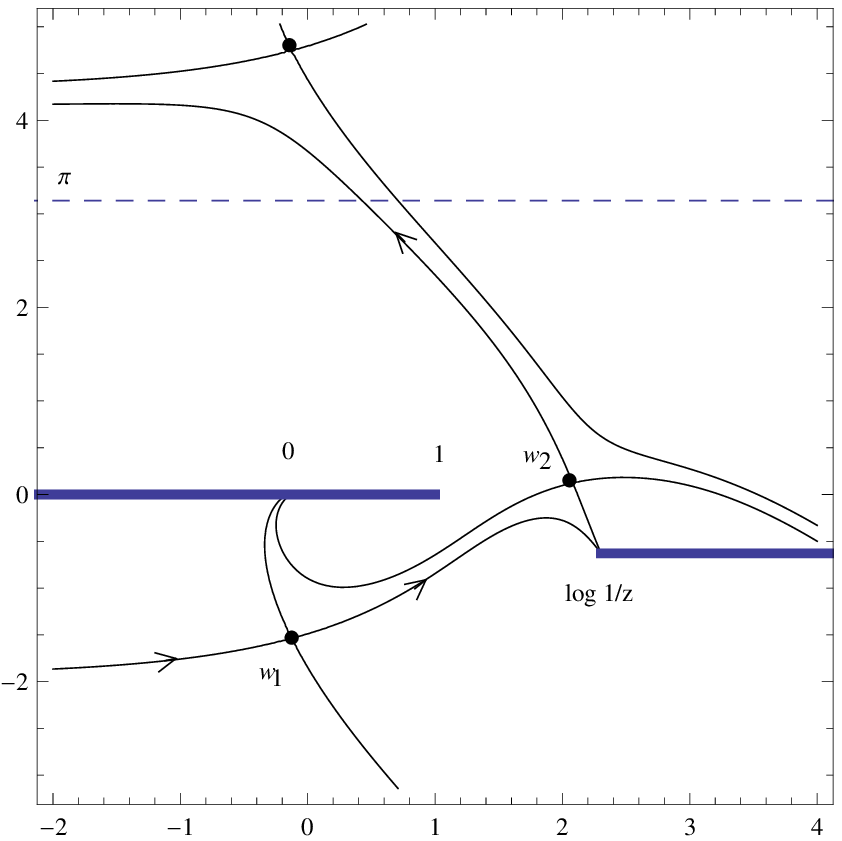}
\vspace{0.3cm}

{\tiny($c$)}\ \includegraphics[width=0.33\textwidth]{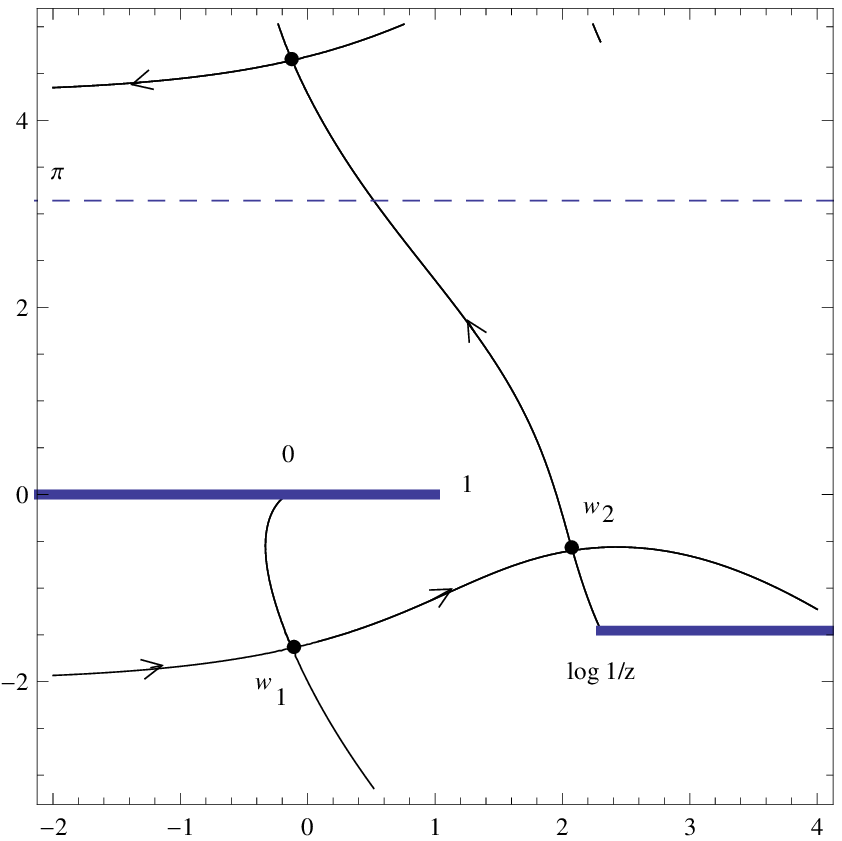}\qquad
{\tiny($d$)}\ \includegraphics[width=0.33\textwidth]{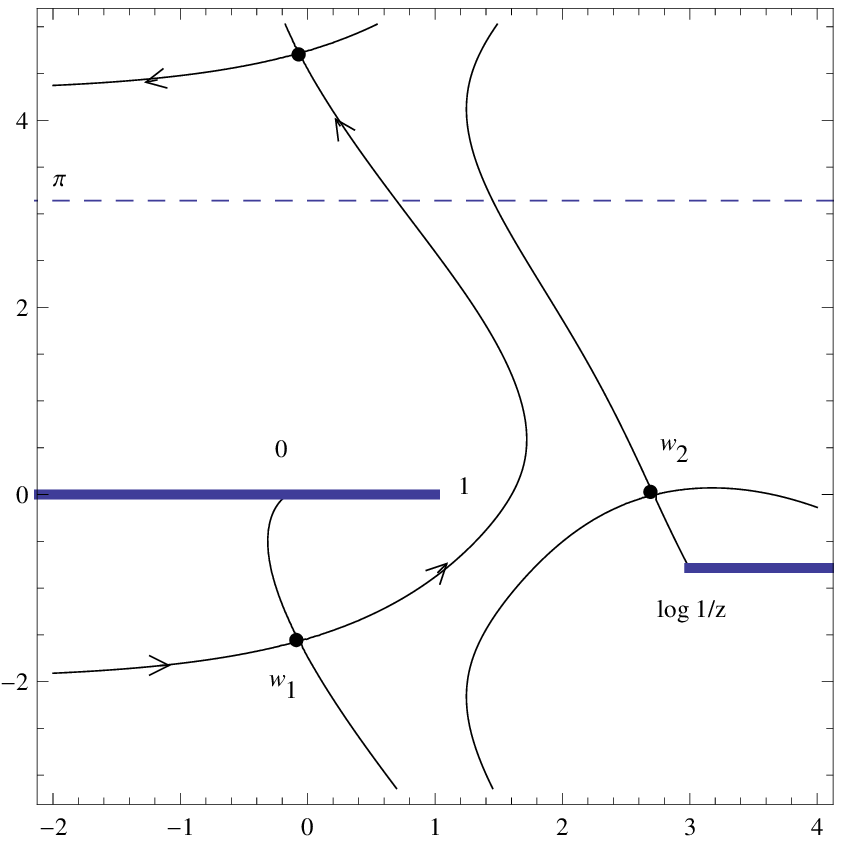}
\caption{\small{Examples of the steepest descent and ascent paths in the $w$-plane for $\alpha=1$ when (a) $|z|=0.30$, $\theta=0.25\pi$, (b) $|z|=0.10$, $\theta=0.20\pi$, (c) $|z|=0.10$, $\theta\doteq 0.46292\pi$ corresponding to a Stokes phenomenon and (d) $|z|=0.10$, $\theta=0.60\pi$. The saddles are denoted by heavy dots and the arrows indicate the integration path. There is a branch cut along $(-\infty,1]$ and from the point $\log\,1/z$ out to infinity.}}
	\end{center}
\end{figure}

Then, from (\ref{e21}), we have the expansion
\bee\label{e26}
F_\alpha(\lambda;z)\sim G_\alpha(\lambda)\left\{\frac{e^{\lambda\psi(t_{s1})} f(t_{s1})}{\sqrt{2\pi \psi''(t_{s1})}}\ \sum_{s=0}^\infty \frac{c_s^{(1)} (\fs)_s}{\lambda^{s+\fr}}+\frac{e^{\lambda\psi(t_{s2})} f(t_{s2})}{\sqrt{2\pi \psi''(t_{s2})}}\ \sum_{s=0}^\infty \frac{c_s^{(2)} (\fs)_s}{\lambda^{s+\fr}}
\right\}
\ee
as $\lambda\ra+\infty$ valid for complex $z$ in $|\arg (1-z)|<\pi$. This expansion will break down in the neighbourhood
of the double saddle points at $z=z_d^\pm$, and also ceases to be valid in a zone 
surrounding $z=0$; see below. The orientation factors $\gamma_1=\gamma_2=0$ and we note that
\[e^{\lambda\psi(t_{sj})}=\bl\{\frac{t_{sj}(1-zt_{sj})}{(t_{sj}-1)^\beta}\br\}^\lambda.\]
If required, an expansion for $G_\alpha(\lambda)$ in inverse powers of $\lambda$ is given in the appendix.

The boundary in the upper-half $z$-plane on which $\Re \psi(t_{s1})=\Re \psi(t_{s2})$ is shown in Fig.~3 for the particular case of $\alpha=1$. This curve has its endpoints $A$ and $B$ at the double saddle points $z_d^\pm$ in (\ref{ezd}).
Below this curve, and also in $\Im (z)<0$, the contribution from the saddle $t_{s1}$ dominates that from the saddle $t_{s2}$; above this curve the saddle $t_{s2}$ is dominant. In the neighbourhood of the curve both contributions need to be taken into account. 

The dashed closed curve surrounding $z=0$ (the enclosed domain is denoted by $\cal D$) and terminating at $z_d^-$ shows the curve on which $\Im \psi(t_{s1})=\Im \psi(t_{s2})$, where a Stokes phenomenon occurs. As one crosses this loop and passes into its interior the saddle $t_{s2}$ disconnects from the integration path to leave only the contribution from the saddle $t_{s1}$; see Fig.~1(d) and Fig.~2(c), (d). 
As a consequence, the expansion of $F_\alpha(\lambda;z)$ inside this loop is given by
\bee\label{e26a}
F_\alpha(\lambda;z)\sim G_\alpha(\lambda)\frac{e^{\lambda\psi(t_{s1})} f(t_{s1})}{\sqrt{2\pi \psi''(t_{s1})}}\ \sum_{s=0}^\infty \frac{c_s^{(1)} (\fs)_s}{\lambda^{s+\fr}}\qquad (z\in \cal D)
\ee
as $\lambda\ra+\infty$. 

In the lower-half $z$-plane, the saddle $t_{s2}$ can pass over onto an adjacent Riemann sheet (in the $t$-plane) and great care must be taken to ensure that one uses continuous branches for the functions $\log\,t_{s2}$ and $\log\,(t_{s2}-1)$. The curves on which $t_{s2}$ and $t_{s2}-1$ pass onto adjacent sheets are indicated in Fig.~3 by the dotted curves issuing from $z=0$ and $z=1$, respectively. To the right (resp. left) of the curve issuing from $z=0$ (resp. $z=1$), $t_{s2}$ (resp. $t_{s2}-1$) lies on the principal sheet. It must be emphasised that all the curves in Fig.~3 depend on the value of the parameter $\alpha>0$. 
\begin{figure}[t]
	\begin{center}
\includegraphics[width=0.55\textwidth]{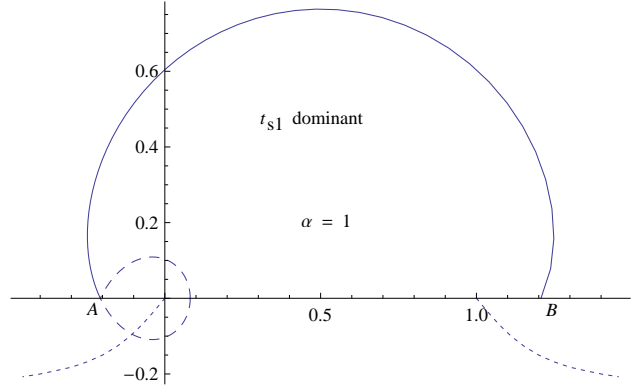}
\caption{\small{The curve (shown solid) in the $z$-plane on which $\Re \psi(t_{s1})=\Re \psi(t_{s2})$ when $\alpha=1$. The endpoints $A$ and $B$ are situated at $z_d^-\doteq -0.2071$ and $z_d^+\doteq 1.0271$ corresponding to the double saddle points in the $t$-plane. The dashed loop surrounding $z=0$ shows the curve on which 
$\Im \psi(t_{s1})=\Im \psi(t_{s2})$, where
a Stokes phenomenon occurs; the interior of this domain is labelled ${\cal D}$. The dotted curves issuing from $z=0$ and $z=1$ indicate where  $t_{s2}$ and $t_{s2}-1$, respectively, pass onto an adjacent Riemann sheet in the $t$-plane.}}
	\end{center}
\end{figure}

The results of numerical computations carried out with {\it Mathematica} are presented in Tables \ref{t21} and \ref{t22}. Table \ref{t21} shows the absolute relative error\footnote{We have adopted the convention in the tables of writing $x(y)$ for $x\times 10^y$.} in the computation of $F_\alpha(\lambda;z)$ as a function of the truncation index $s$ for $\alpha=1$, $\lambda=80$ and $z=0.50e^{i\theta}$ using the expansion (\ref{e26}). 
We note that when $\theta=0$ the value of $F_\alpha(\lambda;\fs)$ is given by (\ref{e14}). 
The coefficients $c_s^{(1,2)}$ were derived from (\ref{e25w}) and the high-precision evaluation of $F_\alpha(\lambda;z)$ obtained by the routine Hypergeometric2F1 in {\it Mathematica}. Table \ref{t22} shows the absolute relative error as a function of $\theta$ for different $|z|$ with the same values of $\alpha$ and $\lambda$ and truncation index $s=2$. In this case the coefficients $c_s^{(1,2)}$ can be obtained alternatively
from (\ref{e25c}). In the column corresponding to $|z|=0.06$, all the values of $z$ lie in the domain ${\cal D}$ in which only the saddle $t_{s1}$ contributes to $F_\alpha(\lambda;z)$. In the column corresponding to $|z|=0.25$,
the error is seen to become relatively large when $\theta=\pi$. This is due to the fact that $z=-0.25$ lies close to the value $z_d^-=-0.2071$, which corresponds to the formation of a double saddle point where the expansion (\ref{e26}) ceases to be valid.

\begin{table}[th]
\caption{\footnotesize{Values of the absolute relative error in the computation of $F_\alpha(\lambda;z)$ for different truncation index $s$ in the expansions (\ref{e26}) and (\ref{e26a}) when $\lambda=80$, $\alpha=1$ and $z=0.50e^{i\theta}$.}}\label{t21}
\begin{center}
\begin{tabular}{|l||l|l|l|l|l|}
\hline
&&&&&\\[-0.25cm]
\mcol{1}{|c||}{$s$} & \mcol{1}{c|}{$\theta=0$} & \mcol{1}{c|}{$\theta=0.25\pi$} & \mcol{1}{c|}{$\theta=0.50\pi$}
& \mcol{1}{c|}{$\theta=0.75\pi$} & \mcol{1}{c|}{$\theta=\pi$}\\
\hline
&&&&&\\[-0.3cm]
0 & $1.042(-08)$ & $1.865(-03)$ & $3.251(-03)$ & $9.873(-04)$ & $6.525(-04)$\\
1 & $5.158(-07)$ & $2.927(-06)$ & $1.619(-05)$ & $1.219(-05)$ & $3.182(-05)$\\
2 & $2.731(-08)$ & $3.390(-08)$ & $3.139(-07)$ & $4.098(-07)$ & $1.617(-06)$\\
3 & $3.050(-11)$ & $7.929(-10)$ & $7.240(-09)$ & $1.910(-08)$ & $1.213(-07)$\\
4 & $2.166(-12)$ & $9.013(-12)$ & $2.263(-10)$ & $1.175(-09)$ & $1.195(-08)$\\
5 & $1.473(-15)$ & $2.314(-13)$ & $8.728(-12)$ & $8.961(-11)$ & $1.464(-09)$\\
[.2cm]\hline\end{tabular}
\end{center}
\end{table}
\begin{table}[th]
\caption{\footnotesize{Values of the absolute relative error in the computation of $F_\alpha(\lambda;z)$ for different $\theta$ and $|z|$ in the expansions (\ref{e26}) and (\ref{e26a}) when $\lambda=80$, $\alpha=1$ and truncation index $s=2$.}}\label{t22}
\begin{center}
\begin{tabular}{|r||l|l|l|l|c|c|}
\hline
&&&&&&\\[-0.25cm]
\mcol{1}{|c||}{$\theta/\pi$} & \mcol{1}{c|}{$|z|=0.06$} & \mcol{1}{c|}{$|z|=0.25$} & \mcol{1}{c|}{$|z|=0.50$}
& \mcol{1}{c|}{$|z|=0.75$} & \mcol{1}{c|}{$|z|=1.00$} & \mcol{1}{c|}{$|z|=1.20$}\\
\hline
&&&&&&\\[-0.3cm]
0       & $1.340(-08)$ & $3.566(-08)$ & $2.731(-08)$ & $3.566(-08)$ & $-\!\!-$ & $-\!\!-$\\
0.25    & $4.546(-08)$ & $7.760(-08)$ & $3.390(-08)$ & $4.584(-08)$ & $9.869(-08)$ & $1.265(-08)$\\
0.50    & $2.064(-07)$ & $5.200(-07)$ & $3.139(-07)$ & $2.347(-08)$ & $1.150(-08)$ & $7.116(-09)$\\
0.75    & $8.228(-07)$ & $1.159(-05)$ & $4.098(-07)$ & $9.420(-08)$ & $3.681(-08)$ & $2.153(-08)$\\
1.00    & $1.522(-06)$ & $7.823(-03)$ & $1.617(-06)$ & $1.860(-07)$ & $5.337(-08)$ & $2.547(-08)$\\
$-0.25$ & $4.848(-08)$ & $4.754(-08)$ & $1.300(-08)$ & $4.925(-09)$ & $1.171(-08)$ & $1.265(-08)$\\
$-0.50$ & $1.597(-07)$ & $2.028(-07)$ & $5.965(-08)$ & $2.347(-08)$ & $1.150(-08)$ & $7.116(-09)$\\
$-0.75$ & $6.302(-07)$ & $3.990(-06)$ & $4.098(-07)$ & $9.420(-08)$ & $4.681(-08)$ & $2.153(-08)$\\
[.2cm]\hline\end{tabular}
\end{center}
\end{table}
\vspace{0.4cm}

\noindent{\bf 2.3\ \ The expansion in the general case}\ \ \ The hypergeometric function in (\ref{e11}) has the integral representation from (\ref{e2int}) given by 
\[F\left(\begin{array}{c}a-\lambda, b+\lambda\\c+i\alpha\lambda\end{array}\!;z\right)=
\frac{\g(c+i\alpha\lambda) \g(1+b-c+\lambda\beta)}{2\pi i\,\g(b+\lambda)} \int_0^{(1+)} f(t) e^{\lambda\psi(t)}dt,\]
where the amplitude function $f(t)$ is now given by
\bee\label{e230}
f(t)=\frac{t^{b-1}(t-1)^{c-b-1}}{(1-zt)^a}~.
\ee
The phase function $\psi(t)$ is as in (\ref{e22}) and consequently the distribution of the saddle points remains the same. It therefore follows that the expansion of $F(a-\lambda,b+\lambda;c+i\alpha\lambda;z)$ for $\lambda\ra+\infty$ is given by (\ref{e26}) and (\ref{e26a}) with $f(t)$ replaced by (\ref{e230}) and
the coefficients $c_s^{(1,2)}$ determined from either (\ref{e25c}) or (\ref{e25w}).
The function $G_\alpha(\lambda)$ is replaced by $\g(c+i\alpha\lambda)\g(1+b-c+\lambda\beta)/\g(b+\lambda)$.

\vspace{0.6cm}

\begin{center}
{\bf 3. \  The expansion when $z=z_d^-$ for $\lambda\ra+\infty$}
\end{center}
\setcounter{section}{3}
\setcounter{equation}{0}
\renewcommand{\theequation}{\arabic{section}.\arabic{equation}}
When $z=z_d^-=\fs(1-\sqrt{1+\alpha^2})$, it is seen from (\ref{e24}) that the saddles $t_{s1}$ and $t_{s2}$ coalesce to form a double saddle at the point
\bee\label{e30}
t_d=\frac{z_d^-+\fs i\alpha}{(1+i\alpha)z_d^-}=\frac{1-\sqrt{1+\alpha^2}+i\alpha}{(1+i\alpha)(1-\sqrt{1+\alpha^2})}.
\ee
In the neighbourhood of the point $z=z_d^-$ the expansions in (\ref{e26}) and (\ref{e26a}) break down. In this section we determine the expansion of $F_\alpha(\lambda;z)$ and also that of the general case in (\ref{e11}) valid at the coalescence point $z=z_d^-$.  The integration path when $z=z_d^-$ is typically as shown in Fig.~4.
\begin{figure}[t]
	\begin{center}
{\tiny($a$)}\ \includegraphics[width=0.35\textwidth]{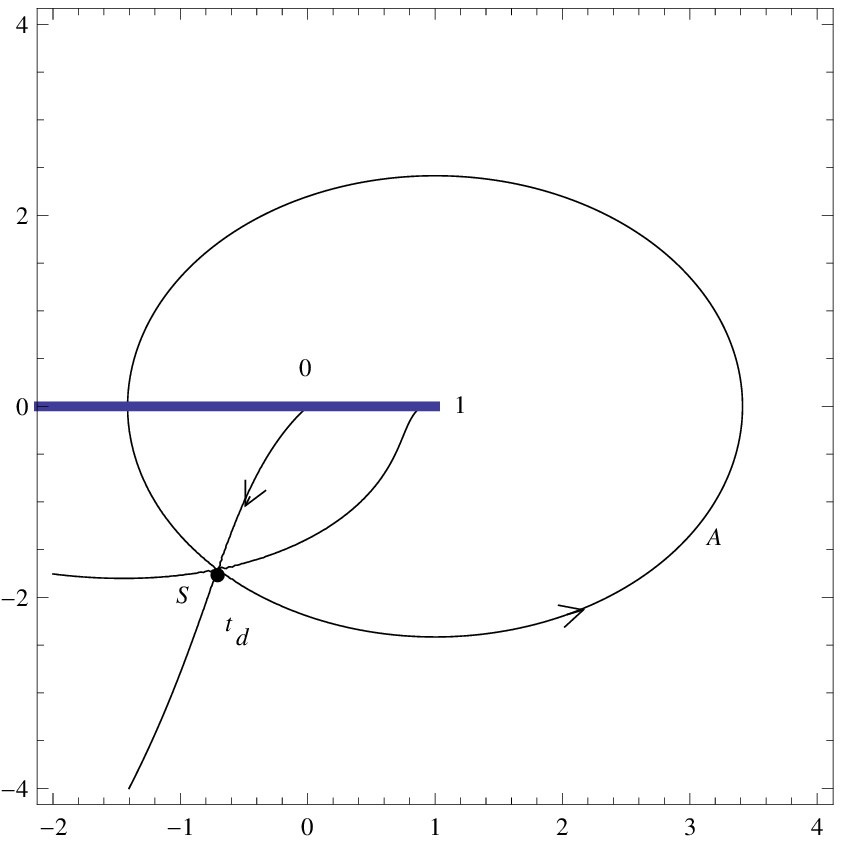}\qquad
{\tiny($b$)}\ \includegraphics[width=0.35\textwidth]{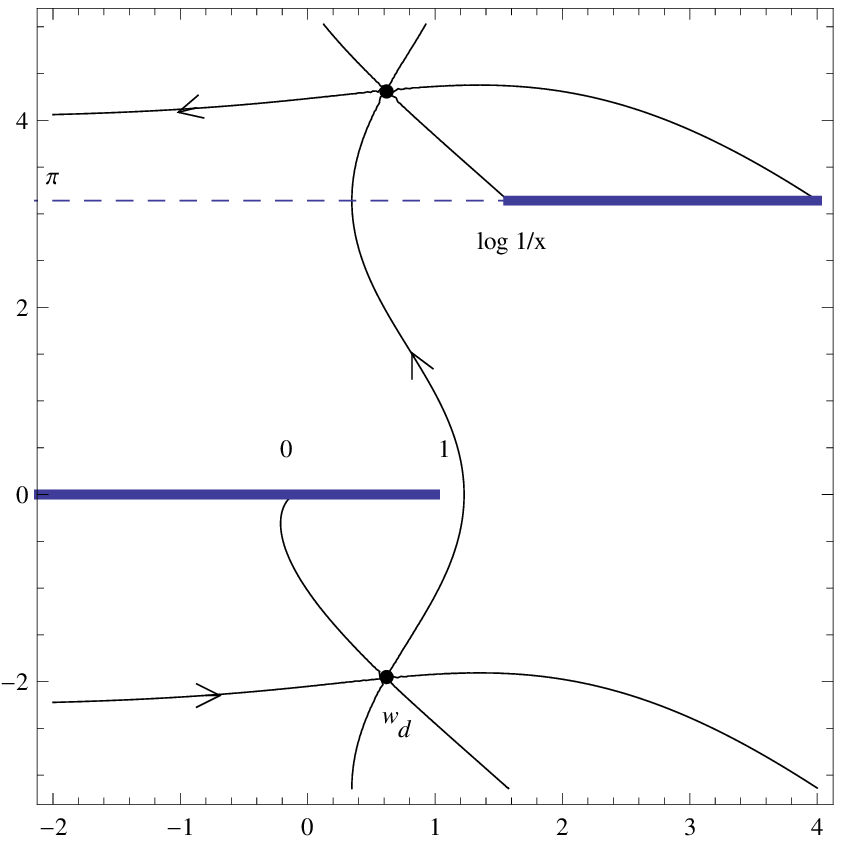}
\caption{\small{An example of the steepest descent and ascent paths when $z=\frac{1}{2}(1-\sqrt{1+\alpha^2})$ and $\alpha=1$. (a) In the $t$-plane with the double saddle at $t_d$. The steepest descent path $SA$ passes over onto the adjacent Riemann sheet and spirals into the origin. (b) The same situation viewed in the $w$-plane where $w_d=\log\,t_d$. The saddles are denoted by heavy dots and the arrows indicate the integration path. There is a branch cut along $(-\infty,1]$ and from the point $1/z$ out to infinity (not shown in (a)).}}
	\end{center}
\end{figure}
\vspace{0.4cm}

\noindent{\bf 3.1\ \ The expansion of $F_\alpha(\lambda;z_d^-)$}.\ \ \ If we put
\[-u=\psi(t)-\psi(t_d)=A\tau^3+B\tau^4+C\tau^5+ D\tau^6+\ldots ,\qquad \tau:=t-t_d\]
we find the coefficients
\[A=\frac{-(1+i\alpha)^3}{6\alpha^2\sqrt{1+\alpha^2}}\bl\{\alpha^2(-3+\sqrt{1+\alpha^2})+4(-1+\sqrt{1+\alpha^2})\br\},\]
\[B=\frac{(1+i\alpha)^3(\alpha+2i)}{8\alpha^3(1-i\alpha)}\bl\{\alpha^4-4\alpha^2(-2+\sqrt{1+\alpha^2})-8(-1+\sqrt{1+\alpha^2})\br\},\]
\[C=\frac{(1+i\alpha)^4(6-5i\alpha-\alpha^2)}{20\alpha^4(1-i\alpha)\sqrt{1+\alpha^2}}\bl\{\alpha^4(-5+
\sqrt{1+\alpha^2})+4\alpha^2(-5+3\sqrt{1+\alpha^2})+16(-1+\sqrt{1+\alpha^2})\br\},\]
\[D=\frac{-i(1+i\alpha)^4}{24\alpha^5(1-i\alpha)^2}(8-9i\alpha-3\alpha^2)\bl\{\alpha^6-6\alpha^4(-3+\sqrt{1+\alpha^2})-16\alpha^2(-3+2\sqrt{1+\alpha^2})\]
\[\hspace{8cm}-32(-1+\sqrt{1+\alpha^2})\br\}, \ldots~.\]
Inversion yields
\bee\label{e440}
\tau(w)=\frac{w^{1/3}}{A^{1/3}}-\frac{Bw^{2/3}}{9A^{5/3}}+\frac{(B^2-AC)w}{3A^3}
-\frac{w^{4/3}}{81A^{13/3}}\,(35B^3-63ABC+27A^2D)+O(w^{5/3}),
\ee
where $w=ue^{\pi i}$ on the path $OS$ and $w=ue^{-\pi i}$ on the path $SA$ in Fig.~4(a).
In addition, with 
\bee\label{e40}
T:=t_d-1=\frac{i\alpha \sqrt{1+\alpha^2}}{(1+i\alpha)(1-\sqrt{1+\alpha^2})},
\ee 
we have 
\[f(t)=\frac{1}{(t-1)}=\frac{1}{T}\bl(1+\frac{\tau}{T}\br)^{-1}=\frac{1}{T}\bl\{1-\frac{w^{1/3}}{A^{1/3}T}+\frac{w^{2/3}}{A^{2/3}T^2}\bl(1+\frac{BT}{3A}\br)\hspace{2cm}\]
\[
\hspace{5cm}-\frac{w}{AT^3}\bl(1+\frac{B^2T^2}{3A^2}+\frac{T}{3A}(2B-CT)\br)+O(w^{4/3})\br\}.
\]

Upon differentiation of $\tau(w)$, we then obtain the expansion
\bee\label{e3fd}
\frac{1}{t-1}\,\frac{d\tau}{dw}=\sum_{m=0}^\infty {\cal B}_m(\alpha) w^{(m-2)/3}
\ee
valid in a neighbourhood of $w=0$ ($t=t_d$), where
\[{\cal B}_0(\alpha)=\frac{1}{3A^{1/3}T},\quad {\cal B}_1(\alpha)=-\frac{1}{9A^{2/3}T^2}\,\bl(3+\frac{2BT}{A}\br),\]
\[{\cal B}_2(\alpha)=\frac{1}{3AT^3}\,\bl\{1+\frac{BT}{A}\bl(1+\frac{BT}{A}\br)-\frac{CT^2}{A}\br\},\]
\bee\label{ecoeffs}
{\cal B}_3(\alpha)=-\frac{1}{243A^{4/3}T^4}\bl\{81+140\frac{B^3T^3}{A^3}+\frac{126BT^2}{A^2}(B-2CT)+\frac{108T}{A}(B-CT+DT^2)\br\}, \ldots~.
\ee
The coefficients $A$, $B$, $C$, $D$ and the quantity $T$ are defined above in terms of the parameter $\alpha$.

Then, from (\ref{e21}) and (\ref{e3fd}), we obtain
\begin{eqnarray*}
F_\alpha(\lambda;z_d^-)&=&\frac{G_\alpha(\lambda) e^{\lambda\psi(t_d)}}{2\pi i} \int_0^\infty e^{-\lambda u}
\bl\{\frac{1}{t-1}\,\left.\frac{dt}{du}\right|_{ue^{-\pi i}}-\frac{1}{t-1}\,\left.\frac{dt}{du}\right|_{ue^{\pi i}}\br\}du\\
&\sim&-\frac{G_\alpha(\lambda)}{\pi}\,e^{\lambda\psi(t_d)}\int_0^\infty e^{-\lambda u}\sum_{m=0}^\infty {\cal B}_m(\alpha) u^{(m-2)/3}\sin \pi(\f{1}{3}m+\f{1}{3})\,du
\end{eqnarray*}
for $\lambda\ra+\infty$.
This therefore produces the expansion
\bee\label{e41}
F_\alpha(\lambda;z_d^-)\sim -\frac{G_\alpha(\lambda)}{\pi}\,e^{\lambda \psi(t_d)} \sum_{m=0}^\infty \frac{{\cal B}_m(\alpha) \g(\f{1}{3}m+\f{1}{3})}{\lambda^{(m+1)/3}}\,\sin \pi(\f{1}{3}m+\f{1}{3})\qquad(\lambda\ra+\infty).
\ee
Since $z_d^-$ is real it immediately follows from (\ref{e21a}) that
\bee\label{e41a}
F_{-\alpha}(\lambda;z_d^-)={\overline F}_\alpha(\lambda;z_d^-).
\ee

Due to the complexity of the coefficients it is not practical to present their explicit dependence on the parameter $\alpha$ for more than the first three terms in the expansion (\ref{e41}). If, however, $\alpha$ is given a numerical value then the inversion process can be carried out with {\it Mathematica} to many more terms.
In the particular case $\alpha=1$, the values of the coefficients ${\cal B}_m(\alpha)$ are tabulated in Table \ref{t31} for $m\leq 10$; we observe that the values of ${\cal B}_2(\alpha)$, ${\cal B}_5(\alpha)$, ${\cal B}_8(\alpha), \ldots$ are not required. Values of $F_\alpha(\lambda;z_d^-)$ and its asymptotic estimate from (\ref{e41}) with $m\leq 10$ (sub-optimal truncation) are given in Table \ref{t32}.
\begin{table}[th]
\caption{\footnotesize{Values of the coefficients ${\cal B}_m(\alpha)$ to 10dp for $m\leq 10$ when $\alpha=1$.}}\label{t31}
\begin{center}
\begin{tabular}{|r|l||r|l|}
\hline
&&&\\[-0.3cm]
\mcol{1}{|c|}{$m$} & \mcol{1}{c||}{${\cal B}_m(\alpha)$} & \mcol{1}{c|}{$m$} & \mcol{1}{c|}{${\cal B}_m(\alpha)$}\\
\hline
&&&\\[-0.3cm]
0 & $+1.1210852199+0.3003938793i$ & 6 &    $+0.0251806251+0.0067471281i$\\
1 & $+0.2166214717+0.8084423383i$ & 7 &    $-0.0018434584-0.0068798804i$\\
3 & $+0.0700678262+0.0187746174i$ & 9 &    $+0.0044976298+0.0012051363i$\\
4 & $+0.0082200112+0.0306774994i$ & 10&    $-0.0001979294-0.0007386825i$\\
[.2cm]\hline\end{tabular}
\end{center}
\end{table}
\begin{table}[th]
\caption{\footnotesize{Values of $F_\alpha(\lambda;z_d^-)$ and the asymptotic approximation (\ref{e41}) with $m\leq 10$ when $\alpha=1$.}}\label{t32}
\begin{center}
\begin{tabular}{|r|l|l|}
\hline
&&\\[-0.35cm]
\mcol{1}{|c|}{$\lambda$} & \mcol{1}{c|}{$F_\alpha(\lambda;z_d^-)$} & \mcol{1}{c|}{Asymptotic\ value}\\
\hline
&&\\[-0.3cm]
10 & $-1.360471986485-1.001859460942i$ & $-1.3604719{\bf 3}7691-1.0018594{\bf 3}5166i$\\
20 & $+0.594081476620+1.768757533401i$ & $+0.59408147{\bf 5}031+1.76875753{\bf 0}253i$\\
40 & $-1.652642124543+1.239971720906i$ & $-1.652642124{\bf 3}73+1.239971720{\bf 7}42i$\\
60 & $-1.790341054904-1.271041237802i$ & $-1.790341054{\bf 8}61-1.271041237{\bf 7}78i$\\
80 & $+0.702398854959-2.183378031108i$ & $+0.70239885495{\bf 5}-2.1833780310{\bf 9}2i$\\
100& $+2.373096383290+0.016497731931i$ & $+2.3730963832{\bf 8}4+0.01649773193{\bf 2}i$\\
[.2cm]\hline\end{tabular}
\end{center}
\end{table}

A uniform approximation for $F_\alpha(\lambda;z)$ for $z\simeq z_d^-$ could be obtained by making the standard cubic transformation (see \cite[(2.4.18)]{DLMF}) to $\psi(t)$ in (\ref{e21}). We do not pursue this further here.
\vspace{0.4cm}

\noindent{\bf 3.2\ \ The expansion in the general case when $z=z_d^-$}.\ \ \ The expansion of the hypergeometric function in (\ref{e11}) for general $a$, $b$ and $c$ when $z=z_d^-$ follows a similar procedure to that in the specific case of $F_\alpha(\lambda;z)$. The amplitude function $f(t)$ is now given by (\ref{e230}), 
which may be expressed in the neighbourhood of the double saddle $t_d$ as (with $t=t_d+\tau$)
\[f(t)=\frac{t^{b-1}(t-1)^{c-b-1}}{(1-z_d^-t)^a}=\frac{1}{{\hat T}}\,(1+\frac{\tau}{T})^{c-b-1}\frac{(1+\frac{\tau}{t_d})^{b-1}}{(1-\frac{\kappa\tau}{t_d})^a},\]
where $T$ is defined in (\ref{e40}) and
\[\kappa:=\frac{z_d^- t_d}{1-z_d^- t_d}=\frac{1-\sqrt{1+\alpha^2}+i\alpha}{1+\sqrt{1+\alpha^2}+i\alpha},
\qquad {\hat T}^{-1}:= \frac{t_d^{b-1} (t_d-1)^{c-b-1}}{(1-z_d^- t_d)^a}~.\]
Expansion about $\tau=0$, followed by use of (\ref{e440}) to express $\tau$ in terms of $w$, produces
\[\log\,f(t)=-\log\,{\hat T}+\sum_{n=1}^\infty \gamma_n \bl(\frac{w}{A}\br)^{\!\!n/3},\]
where
\[\gamma_1=\frac{c-b-1}{T}+\frac{b-1+\kappa a}{t_d},\qquad \gamma_2=-\frac{b-1-\kappa^2a}{2t_d^2}-\frac{c-b-1}{2T^2}-\frac{\gamma_1 B}{3A},\]
\[\gamma_3=\frac{(B^2-3AC)\gamma_1}{9A^2}-\frac{2B\gamma_2}{3A}+\frac{b-1+\kappa^3a}{3t_d^3}+\frac{c-b-1}{3T^3},\, \ldots\, . \]
Application of Lemma 1 in the appendix therefore shows that
\[f(t)=\frac{1}{{\hat T}}\bl\{1+\sum_{n=0}^\infty D_n \bl(\frac{w}{A}\br)^{\!\!n/3}\br\},\]
where
\[D_1=\gamma_1,\quad D_2=\fs\gamma_1^2+\gamma_2,\quad D_3=\f{1}{6}\gamma_1^3+\gamma_1\gamma_2+\gamma_3, \ldots\,.\]

Then, from (\ref{e440}) we obtain the result
\[f(t) \frac{dt}{dw}=\sum_{m=0}^\infty {\hat B}_m(\alpha) w^{(m-2)/3}\]
valid near $w=0$ ($t=t_d$), where the first three contributory coefficients are\footnote{We omit the coefficient ${\hat B}_2(\alpha)$ as this is not required.}
\[{\hat B}_0(\alpha)=\frac{1}{3A^{1/3} {\hat T}},\qquad {\hat B}_1(\alpha)=\frac{1}{9A^{2/3} {\hat T}}\bl\{3\bl(\frac{c-b-1}{T}+\frac{b-1+\kappa a}{t_d}\br)-\frac{2B}{A}\br\},\]
\bee\label{b2}
{\hat B}_3(\alpha)=-\frac{1}{243A^{4/3}{\hat T}}\bl\{\frac{4}{A^3}(35B^3-63ABC+27A^2D)-\frac{81}{A^2}(B^2-AC)D_1+\frac{54BD_2}{A}-81D_3\br\}.
\ee
From this we find the expansion
\[F\left(\begin{array}{c}a-\lambda, b+\lambda\\c+i\alpha\lambda\end{array}\!;z\right)\sim
\frac{\g(c+i\alpha\lambda) \g(1+b-c+\lambda\beta)}{\g(b+\lambda)}\,\frac{e^{\lambda\psi(t_d)}}{\pi}\hspace{4cm}\] 
\bee\label{e441}
\hspace{3cm}\times \sum_{m=0}^\infty {\hat B}_m(\alpha) \frac{\g(\f{1}{3}m+\f{1}{3})}{\lambda^{(m+1)/3}}\,\sin \pi(\f{1}{3}m+\f{1}{3})\qquad(\lambda\ra+\infty).
\ee
The complexity of the higher contributory coefficients ${\hat B}_m(\alpha)$ ($m\geq 4$) is such that there appears to be little value in their presentation, although in specific cases with numerical values for $a$, $b$, $c$ and $\alpha$ it would be quite feasible to continue the inversion process to higher order. 
It can be verified with some effort that when $a=0$, $b=c=1$ we have ${\hat T}=T$ and the coefficients in (\ref{b2}) reduce to those given in (\ref{ecoeffs}).

\vspace{0.6cm}

\begin{center}
{\bf 4. \  Application to the expansion of the Legendre functions}
\end{center}
\setcounter{section}{4}
\setcounter{equation}{0}
\renewcommand{\theequation}{\arabic{section}.\arabic{equation}}
From (\ref{e12}) and (\ref{e13}) we have the Legendre functions of degree $\lambda$ and order $-i\alpha\lambda$,
where $\alpha>0$, $\lambda>0$, given by
\bee\label{ee41}
P_\lambda^{-i\alpha\lambda}(x)=\frac{1}{\g(1+i\alpha\lambda)}\,\bl(\frac{x-1}{x+1}\br)^{i\alpha\lambda/2}\,F_\alpha(\lambda;\fs-\fs x),
\ee

\[e^{-\pi\alpha\lambda}Q_\lambda^{-i\alpha\lambda}(x)=\frac{\g(-i\alpha\lambda)}{2}\,\bl(\frac{x-1}{x+1}\br)^{i\alpha\lambda/2} F_\alpha(\lambda;\fs-\fs x)\hspace{4cm}\]
\bee\label{ee42}
\hspace{4cm}+\frac{\g(i\alpha\lambda)\g(1+\lambda\beta)}{2\g(1+\lambda{\overline \beta})} \,\bl(\frac{x-1}{x+1}\br)^{-i\alpha\lambda/2} F_{-\alpha}(\lambda;\fs-\fs x),
\ee
where $F_\alpha(\lambda;z)$ is defined in (\ref{e20}) and $F_{-\alpha}(\lambda;\fs-\fs x)$ is given by the conjugate relation (\ref{e21a}). These functions are defined\footnote{In {\it Mathematica} they are obtained numerically by use of the `type-3' Legendre functions.} in the complex $x$-plane cut along $(-\infty,1]$. 
The expansions for $P_\lambda^{\pm i\alpha\lambda}(x)$ and $Q_\lambda^{\pm i\alpha\lambda}(x)$ then follow from that of $F_\alpha(\lambda;z)$ given in (\ref{e26}) and (\ref{e26a}).

In the special case $x=\sqrt{1+\alpha^2}$, the argument of $F_{\pm\alpha}(\lambda;\fs-\fs x)$ is equal to 
\[z_d^-=\fs(1-\sqrt{1+\alpha^2})\]
in (\ref{ezd}). This corresponds to the coincidence of the two saddle points associated with the integral for $F_\alpha(\lambda;z)$. From the expansion in (\ref{e41}), it therefore follows that
\bee\label{ee43}
P_\lambda^{-i\alpha\lambda}(\sqrt{1+\alpha^2})\sim -\frac{\g(1+\lambda\beta)}{\g(1+\lambda)}\,\bl(\frac{\sqrt{1+\alpha^2}-1}{\sqrt{1+\alpha^2}+1}\br)^{i\alpha\lambda/2} \frac{e^{\lambda\psi(t_d)}}{\pi}\,S(\lambda;\alpha)
\ee 
and
\bee\label{ee44}
Q_\lambda^{-i\alpha\lambda}(\sqrt{1+\alpha^2})\sim \frac{e^{\pi\alpha\lambda}}{\sinh \pi\alpha\lambda}\, \frac{\g(1+\lambda\beta)}{\g(1+\lambda)}\,\Im\bl\{\bl(\frac{\sqrt{1+\alpha^2}-1}{\sqrt{1+\alpha^2}+1}\br)^{i\alpha\lambda/2} e^{\lambda\psi(t_d)}\,S(\lambda;\alpha)\br\}
\ee
as $\lambda\ra+\infty$, where
\[S(\lambda;\alpha):=\sum_{m=0}^\infty \frac{{\cal B}_m(\alpha) \g(\f{1}{3}m+\f{1}{3})}{\lambda^{(m+1)/3}}\,\sin \pi(\f{1}{3}m+\f{1}{3}),\]
$t_d$ is given in (\ref{e30}) and $\psi(t)$ is defined in (\ref{e22}).

The expansions for the functions with positive imaginary order then follow from
\[P_\lambda^{i\alpha\lambda}(\sqrt{1+\alpha^2})={\overline P}_\lambda^{\,-i\alpha\lambda}(\sqrt{1+\alpha^2})\]
and \cite[(14.9.14), (14.3.10)]{DLMF}
\[Q_\lambda^{i\alpha\lambda}(\sqrt{1+\alpha^2})=\frac{\g(1+\lambda{\overline \beta})}{\g(1+\lambda\beta)}\,e^{-2\pi\alpha\lambda}\,Q_\lambda^{-i\alpha\lambda}(\sqrt{1+\alpha^2}).\]

\vspace{0.6cm}

\begin{center}
{\bf Appendix: \ The expansion of $G_\alpha(\lambda)$}
\end{center}
\setcounter{section}{1}
\setcounter{equation}{0}
\renewcommand{\theequation}{\Alph{section}.\arabic{equation}}
In this appendix we consider the expansion of $G_\alpha(\lambda)$ in (\ref{e23}) in inverse powers of $\lambda$. 
This is given for completeness as the main asymptotic problem under consideration is the large-$\lambda$ expansion of the integral appearing in (\ref{e21}). 

It is sufficient to consider $\alpha> 0$ since the value of $G_\alpha(\lambda)$ for $\alpha<0$ is given by its conjugate (when $\lambda>0$).
We use the well-known expansion for $\log\,\g(z)$ as $|z|\ra\infty$ is \cite[p.~141]{DLMF}
\[\log\,\g(z)\sim(z-\frac{1}{2}) \log\,z-z+\frac{1}{2}\log\,2\pi+\sum_{k=1}^\infty \frac{B_{2k}}{2k(2k-1)z^{2k-1}}\qquad(|\arg\,z|<\pi),\]
where $B_{2k}$ denote the even-order Bernoulli numbers. In addition, we have the following lemma \cite{CEV}:
\newtheorem{lemma}{Lemma}
\begin{lemma}$\!\!\!.$\ \ Let $S(x)=\sum_{n=1}^\infty a_n x^{-n}$ as $x\ra\infty$ be a given expansion. Then the composition $\exp [S(x)]$ has the asymptotic expansion of the following form\footnote{The coefficients $b_n$ can also be expressed in terms of the complete Bell polynomial $B_n$ in the form $b_n=B_n(1! a_1,2! a_2, \ldots , n! a_n)/n!$.}
\[\exp [S(x)] \sim \sum_{n=0}^\infty b_n x^{-n} \qquad (x\ra\infty)\]
where
\[b_0=1,\qquad b_n=\frac{1}{n}\sum_{k=1}^n ka_k b_{n-k} \qquad(n\geq 1).\]
\end{lemma} 

Then it follows that
\[\log \bl(\frac{G_\alpha(\lambda)}{i\alpha\beta\lambda}\br)=\log \frac{\g(i\alpha\lambda) \g(\lambda\beta)}{G_\alpha(\lambda)}\hspace{8cm}\]
\[\hspace{3cm}\sim \frac{1}{2} \log \frac{2\pi}{\lambda}+(i\alpha\lambda-\frac{1}{2}) \log\i\alpha+(\lambda\beta-\frac{1}{2}) \log (1-i\alpha)+\sum_{k=1}^\infty \frac{A_k}{\lambda^{2k-1}}\]
as $\lambda\ra+\infty$, where
\[A_k:=\frac{B_{2k}}{2k(2k-1)}\bl\{\frac{1}{(i\alpha)^{2k-1}}+\frac{1}{(1-i\alpha)^{2k-1}}-1\br\}.\]
Some straightforward algebra and application of Lemma 1 then produces
\begin{eqnarray}
G_\alpha(\lambda)&\sim& \sqrt{2\pi\lambda\alpha}\,(1+\alpha^2)^{\fr\lambda+\frac{1}{4}} e^{-\lambda\alpha(\fr\pi+\phi)+i\Phi}
\exp\bl\{\sum_{k=1}^\infty\frac{A_k}{\lambda^{2k-1}}\br\}\nonumber\\
&\sim&\sqrt{2\pi\lambda\alpha}\,(1+\alpha^2)^{\fr\lambda+\frac{1}{4}} e^{-\lambda\alpha(\fr\pi+\phi)+i\Phi} \, \bl\{1+\sum_{k=1}^\infty \frac{E_k}{\lambda^k}\br\}\qquad (\lambda\ra+\infty)\label{e27}
\end{eqnarray}
where 
\[\Phi:=\lambda\alpha \log \bl(\frac{\alpha}{\sqrt{1+\alpha^2}}\br)-(\lambda+\fs)\phi+\f{1}{4}\pi\]
with $\phi$ defined in (\ref{e200}).
The first few coefficients $E_k$ are
\[E_1=A_1,\quad E_2=\fs A_1^2,\quad E_3=\f{1}{6}A_1^3+A_2,\quad E_4=\f{1}{24}A_1^4+A_1A_2,\]
\[E_5=\f{1}{120}A_1^5+\fs A_1^2A_2+A_3,\quad 
E_6=\f{1}{720}A_1^6+\f{1}{6}A_1^3A_2+\fs A_2^2+A_1 A_3, \ldots\ .\]


\vspace{0.6cm}

\begin{thebibliography}{99}
\footnotesize{
\bibitem{BH}
N. Bleistein and R.A. Handelsman, Asymptotic Expansions of Integrals, Dover, New York 1986.

\bibitem{CEV}
C.-P. Chen, N. Elezovi\'c and L. Vuk\v{s}i\'c, Asymptotic formulae associated with the Wallis power function and digamma function, J. Classical Anal. {\bf 2} (2013) 151--166.

\bibitem{Cherry}
T.M. Cherry,  Asymptotic expansions for the hypergeometric functions occurring in gas-flow theory, Proc. Roy. Soc. London {\bf A202} (1950) 507--522.

\bibitem{Cvit}
M. Cvitkovi\'c, A.-S. Smith and J. Pande, General asymptotic expansions of the hypergeometric function with two large parameters. arXiv:1602.05146 (2016).

\bibitem{D}
R.B. Dingle, {\it Asymptotic Expansions: Their Derivation and Interpretation}, Academic Press, London 1973.

\bibitem{TMD}
T.M. Dunster, Uniform asymptotic solutions of second-order linear differential equations having a double pole with complex exponent and a coalescing turning point, SIAM J. Math. Anal. {\bf 21} (1990) 1594--1618.

\bibitem{DSJ}
D.S. Jones, Asymptotics of the hypergeometric function, Math. Meth. Appl. Sci. {\bf 24} (2001) 369--389.

\bibitem{Lig}
M.J. Lighthill, The hodograph transformation in trans-sonic flow. II Auxiliary theorems on the hypergeometric functions $\psi_n(\tau)$. Proc. Roy. Soc. London {\bf A191} (1947) 341--351.

\bibitem{OD}
A.B. Olde Daalhuis, Uniform asymptotic expansions for hypergeometric functions. I, II. Anal. Appl. (Singap.) {\bf 1} (2003) 111--120, 121--128.

\bibitem{O}
F.W.J. Olver, {\it Asymptotics and Special Functions}, Academic Press, New York 1974. Reprinted A.K. Peters, Massachussets 1997.

\bibitem{DLMF}
F.W.J. Olver, D.W. Lozier, R.F. Boisvert and C.W. Clark (eds.),    
{\it NIST Handbook of Mathematical Functions}, Cambridge University Press, Cambridge 2010.

\bibitem{PBook}
R.B. Paris, {\it Hadamard Expansions and Hyperasymptotic Evaluation: An Extension of the Method of Steepest Descents}, Cambridge University Press, Cambridge 2011.

\bibitem{P1}
R.B. Paris Asymptotics of the Gauss hypergeometric function with large parameters, I., J. Classical Anal. {\bf 2} (2013) 183--203.

\bibitem{P2}
R.B. Paris Asymptotics of the Gauss hypergeometric function with large parameters, II., J. Classical Anal. {\bf 3} (2013) 1--15.

\bibitem{T}
N.M. Temme, {\it Asymptotic Methods for Integrals}, Series in Analysis vol. 6, World Scientific, New Jersey 2015.

\bibitem{W}
G.N. Watson, Asymptotic expansions of hypergeometric functions, Trans. Cambridge Philos. Soc. {\bf 22} (1918) 277--308.

\bibitem{VSWP}
A. Vainchtein, Y. Starosvetsky, J.D. Wright and R. Perline, Solitary waves in diatomic chains, Phys. Rev. E{\bf 93} (2016) 042210.

\bibitem{Woj}
J. Wojdylo, On the coefficients that arise from Laplace's method, J. Comput. Appl. Math. {\bf 196} (2006) 241--266.

\bibitem{Wr}
J.D. Wright, Private communication (2015).
}
\end{thebibliography}
\end{document}